\documentclass[reqno,12pt,a4paper]{amsart}

\usepackage{mathrsfs}
\usepackage{amssymb}
\usepackage{amsfonts}
\usepackage{latexsym}
\usepackage{amsthm}

\usepackage{graphicx}
\def\lb{\label}

\newcommand{\er}[1]{\textrm{(\ref{#1})}}

\begin{document}


\renewcommand{\theequation}{\arabic{section}.\arabic{equation}}
\theoremstyle{plain}
\newtheorem{theorem}{\bf Theorem}[section]
\newtheorem{lemma}[theorem]{\bf Lemma}
\newtheorem{corollary}[theorem]{\bf Corollary}
\newtheorem{proposition}[theorem]{\bf Proposition}
\newtheorem{definition}[theorem]{\bf Definition}
\newtheorem{remark}[theorem]{\it Remark}

\def\a{\alpha}  \def\cA{{\mathcal A}}     \def\bA{{\bf A}}  \def\mA{{\mathscr A}}
\def\b{\beta}   \def\cB{{\mathcal B}}     \def\bB{{\bf B}}  \def\mB{{\mathscr B}}
\def\g{\gamma}  \def\cC{{\mathcal C}}     \def\bC{{\bf C}}  \def\mC{{\mathscr C}}
\def\G{\Gamma}  \def\cD{{\mathcal D}}     \def\bD{{\bf D}}  \def\mD{{\mathscr D}}
\def\d{\delta}  \def\cE{{\mathcal E}}     \def\bE{{\bf E}}  \def\mE{{\mathscr E}}
\def\D{\Delta}  \def\cF{{\mathcal F}}     \def\bF{{\bf F}}  \def\mF{{\mathscr F}}
\def\c{\chi}    \def\cG{{\mathcal G}}     \def\bG{{\bf G}}  \def\mG{{\mathscr G}}
\def\z{\zeta}   \def\cH{{\mathcal H}}     \def\bH{{\bf H}}  \def\mH{{\mathscr H}}
\def\e{\eta}    \def\cI{{\mathcal I}}     \def\bI{{\bf I}}  \def\mI{{\mathscr I}}
\def\p{\psi}    \def\cJ{{\mathcal J}}     \def\bJ{{\bf J}}  \def\mJ{{\mathscr J}}
\def\vT{\Theta} \def\cK{{\mathcal K}}     \def\bK{{\bf K}}  \def\mK{{\mathscr K}}
\def\k{\kappa}  \def\cL{{\mathcal L}}     \def\bL{{\bf L}}  \def\mL{{\mathscr L}}
\def\l{\lambda} \def\cM{{\mathcal M}}     \def\bM{{\bf M}}  \def\mM{{\mathscr M}}
\def\L{\Lambda} \def\cN{{\mathcal N}}     \def\bN{{\bf N}}  \def\mN{{\mathscr N}}
\def\m{\mu}     \def\cO{{\mathcal O}}     \def\bO{{\bf O}}  \def\mO{{\mathscr O}}
\def\n{\nu}     \def\cP{{\mathcal P}}     \def\bP{{\bf P}}  \def\mP{{\mathscr P}}
\def\r{\rho}    \def\cQ{{\mathcal Q}}     \def\bQ{{\bf Q}}  \def\mQ{{\mathscr Q}}
\def\s{\sigma}  \def\cR{{\mathcal R}}     \def\bR{{\bf R}}  \def\mR{{\mathscr R}}
\def\S{\Sigma}  \def\cS{{\mathcal S}}     \def\bS{{\bf S}}  \def\mS{{\mathscr S}}
\def\t{\tau}    \def\cT{{\mathcal T}}     \def\bT{{\bf T}}  \def\mT{{\mathscr T}}
\def\f{\phi}    \def\cU{{\mathcal U}}     \def\bU{{\bf U}}  \def\mU{{\mathscr U}}
\def\F{\Phi}    \def\cV{{\mathcal V}}     \def\bV{{\bf V}}  \def\mV{{\mathscr V}}
\def\P{\Psi}    \def\cW{{\mathcal W}}     \def\bW{{\bf W}}  \def\mW{{\mathscr W}}
\def\o{\omega}  \def\cX{{\mathcal X}}     \def\bX{{\bf X}}  \def\mX{{\mathscr X}}
\def\x{\xi}     \def\cY{{\mathcal Y}}     \def\bY{{\bf Y}}  \def\mY{{\mathscr Y}}
\def\X{\Xi}     \def\cZ{{\mathcal Z}}     \def\bZ{{\bf Z}}  \def\mZ{{\mathscr Z}}
\def\O{\Omega}

\newcommand{\mc}{\mathscr {c}}

\newcommand{\gA}{\mathfrak{A}}          \newcommand{\ga}{\mathfrak{a}}
\newcommand{\gB}{\mathfrak{B}}          \newcommand{\gb}{\mathfrak{b}}
\newcommand{\gC}{\mathfrak{C}}          \newcommand{\gc}{\mathfrak{c}}
\newcommand{\gD}{\mathfrak{D}}          \newcommand{\gd}{\mathfrak{d}}
\newcommand{\gE}{\mathfrak{E}}
\newcommand{\gF}{\mathfrak{F}}           \newcommand{\gf}{\mathfrak{f}}
\newcommand{\gG}{\mathfrak{G}}           
\newcommand{\gH}{\mathfrak{H}}           \newcommand{\gh}{\mathfrak{h}}
\newcommand{\gI}{\mathfrak{I}}           \newcommand{\gi}{\mathfrak{i}}
\newcommand{\gJ}{\mathfrak{J}}           \newcommand{\gj}{\mathfrak{j}}
\newcommand{\gK}{\mathfrak{K}}            \newcommand{\gk}{\mathfrak{k}}
\newcommand{\gL}{\mathfrak{L}}            \newcommand{\gl}{\mathfrak{l}}
\newcommand{\gM}{\mathfrak{M}}            \newcommand{\gm}{\mathfrak{m}}
\newcommand{\gN}{\mathfrak{N}}            \newcommand{\gn}{\mathfrak{n}}
\newcommand{\gO}{\mathfrak{O}}
\newcommand{\gP}{\mathfrak{P}}             \newcommand{\gp}{\mathfrak{p}}
\newcommand{\gQ}{\mathfrak{Q}}             \newcommand{\gq}{\mathfrak{q}}
\newcommand{\gR}{\mathfrak{R}}             \newcommand{\gr}{\mathfrak{r}}
\newcommand{\gS}{\mathfrak{S}}              \newcommand{\gs}{\mathfrak{s}}
\newcommand{\gT}{\mathfrak{T}}             \newcommand{\gt}{\mathfrak{t}}
\newcommand{\gU}{\mathfrak{U}}             \newcommand{\gu}{\mathfrak{u}}
\newcommand{\gV}{\mathfrak{V}}             \newcommand{\gv}{\mathfrak{v}}
\newcommand{\gW}{\mathfrak{W}}             \newcommand{\gw}{\mathfrak{w}}
\newcommand{\gX}{\mathfrak{X}}               \newcommand{\gx}{\mathfrak{x}}
\newcommand{\gY}{\mathfrak{Y}}              \newcommand{\gy}{\mathfrak{y}}
\newcommand{\gZ}{\mathfrak{Z}}             \newcommand{\gz}{\mathfrak{z}}

\def\ve{\varepsilon}   \def\vt{\vartheta}    \def\vp{\varphi}    \def\vk{\varkappa}

\def\A{{\mathbb A}} \def\B{{\mathbb B}} \def\C{{\mathbb C}}
\def\dD{{\mathbb D}} \def\E{{\mathbb E}} \def\dF{{\mathbb F}} \def\dG{{\mathbb G}} \def\H{{\mathbb H}}\def\I{{\mathbb I}} \def\J{{\mathbb J}} \def\K{{\mathbb K}} \def\dL{{\mathbb L}}\def\M{{\mathbb M}} \def\N{{\mathbb N}} \def\dO{{\mathbb O}} \def\dP{{\mathbb P}} \def\R{{\mathbb R}}\def\S{{\mathbb S}} \def\T{{\mathbb T}} \def\U{{\mathbb U}} \def\V{{\mathbb V}}\def\W{{\mathbb W}} \def\X{{\mathbb X}} \def\Y{{\mathbb Y}} \def\Z{{\mathbb Z}}


\def\la{\leftarrow}              \def\ra{\rightarrow}            \def\Ra{\Rightarrow}
\def\ua{\uparrow}                \def\da{\downarrow}
\def\lra{\leftrightarrow}        \def\Lra{\Leftrightarrow}


\def\lt{\biggl}                  \def\rt{\biggr}
\def\ol{\overline}               \def\wt{\widetilde}
\def\no{\noindent}


\let\ge\geqslant                 \let\le\leqslant
\def\lan{\langle}                \def\ran{\rangle}
\def\/{\over}                    \def\iy{\infty}
\def\sm{\setminus}               \def\es{\emptyset}
\def\ss{\subset}                 \def\ts{\times}
\def\pa{\partial}                \def\os{\oplus}
\def\om{\ominus}                 \def\ev{\equiv}
\def\iint{\int\!\!\!\int}        \def\iintt{\mathop{\int\!\!\int\!\!\dots\!\!\int}\limits}
\def\el2{\ell^{\,2}}             \def\1{1\!\!1}
\def\sh{\sharp}
\def\wh{\widehat}
\def\bs{\backslash}
\def\intl{\int\limits}

\def\na{\mathop{\mathrm{\nabla}}\nolimits}
\def\sh{\mathop{\mathrm{sh}}\nolimits}
\def\ch{\mathop{\mathrm{ch}}\nolimits}
\def\where{\mathop{\mathrm{where}}\nolimits}
\def\all{\mathop{\mathrm{all}}\nolimits}
\def\as{\mathop{\mathrm{as}}\nolimits}
\def\Area{\mathop{\mathrm{Area}}\nolimits}
\def\arg{\mathop{\mathrm{arg}}\nolimits}
\def\const{\mathop{\mathrm{const}}\nolimits}
\def\det{\mathop{\mathrm{det}}\nolimits}
\def\diag{\mathop{\mathrm{diag}}\nolimits}
\def\diam{\mathop{\mathrm{diam}}\nolimits}
\def\dim{\mathop{\mathrm{dim}}\nolimits}
\def\dist{\mathop{\mathrm{dist}}\nolimits}
\def\Im{\mathop{\mathrm{Im}}\nolimits}
\def\Iso{\mathop{\mathrm{Iso}}\nolimits}
\def\Ker{\mathop{\mathrm{Ker}}\nolimits}
\def\Lip{\mathop{\mathrm{Lip}}\nolimits}
\def\rank{\mathop{\mathrm{rank}}\limits}
\def\Ran{\mathop{\mathrm{Ran}}\nolimits}
\def\Re{\mathop{\mathrm{Re}}\nolimits}
\def\Res{\mathop{\mathrm{Res}}\nolimits}
\def\res{\mathop{\mathrm{res}}\limits}
\def\sign{\mathop{\mathrm{sign}}\nolimits}
\def\span{\mathop{\mathrm{span}}\nolimits}
\def\supp{\mathop{\mathrm{supp}}\nolimits}
\def\Tr{\mathop{\mathrm{Tr}}\nolimits}
\def\BBox{\hspace{1mm}\vrule height6pt width5.5pt depth0pt \hspace{6pt}}


\newcommand\nh[2]{\widehat{#1}\vphantom{#1}^{(#2)}}
\def\dia{\diamond}

\def\Oplus{\bigoplus\nolimits}



\def\qqq{\qquad}
\def\qq{\quad}
\let\ge\geqslant
\let\le\leqslant
\let\geq\geqslant
\let\leq\leqslant
\newcommand{\ca}{\begin{cases}}
\newcommand{\ac}{\end{cases}}
\newcommand{\ma}{\begin{pmatrix}}
\newcommand{\am}{\end{pmatrix}}
\renewcommand{\[}{\begin{equation}}
\renewcommand{\]}{\end{equation}}
\def\eq{\begin{equation}}
\def\qe{\end{equation}}
\def\[{\begin{equation}}
\def\bu{\bullet}

\title[{Global estimates of resonances for 1D Dirac operators}]
{Global estimates of resonances for 1D Dirac operators}

\date{\today}

\author[Evgeny Korotyaev]{Evgeny L. Korotyaev}
\address{Mathematical Physics Department, Faculty of Physics, Ulianovskaya 2,
St. Petersburg State University, St. Petersburg, 198904,
 and Pushkin Leningrad State University, Russia,
 \ korotyaev@gmail.com,
}


\begin{abstract}
\no We discuss resonances for 1D massless Dirac operators  with
compactly supported potentials on the line. We estimate the sum of
the negative power of all resonances in terms of the norm of the
potential and the diameter of its support.

\end{abstract}

\maketitle

\vskip 0.25cm
\section {Introduction and main results}
\setcounter{equation}{0}

In this paper we plan to determine  global estimates of resonances
in terms of the potential for massless  Dirac operators  $H$ acting
in $L^2(\R )\os L^2(\R )$ and given by
$$
H=-iJ {d\/ dx}+ V,\qq \qqq  J =\ma 1&0\\ 0&-1\am, \ \ \
\ \ V= \ma 0&q\\
\ol q & 0\am.
$$
 Here  $q$ is a complex-valued , integrable function with compact support
 $\supp q\ss [0,\g]$ for some $\g>0$.
It is well known  that the operator $H$ is self-adjoint (see
Theorem 3.2 in \cite{LM03}) and the spectrum of $H$ is purely
absolutely continuous and covers the real line (see  \cite{DEGM82}).

{\bf Below we consider all functions and the resolvent in upper-half
plane $\C_+$ and we will obtain their analytic extensions into the
whole complex plane $\C$.} Note that we can consider all functions
and the resolvent in lower-half plane $\C_-$ and to obtain their
analytic extensions into the whole complex plane $\C$. The Riemann
surface of the resolvent for the Dirac operator consists of two
not-connected sheets $\C$. In the case of the Schr\"odinger operator
 the corresponding Riemann surface is the Riemann surface  of the function
 $\sqrt \l$.

We consider the fundamental solutions $\p^{\pm}$ of the Dirac
equation
\[
\lb{1.2} -iJ f'+ V(x)f=\l f
\]
under the following conditions
\[
\p^{\pm}(x,\l )=e^{\pm i\l x}e_{\pm},\ \ \ \  x>\g;\ \ \  \ \  \ \ \
\ \ \ \vp^{\pm}(x,\l )=e^{\pm i\l x}e_{\pm},\ \ \ \  x<0.
\]
where the vectors $e_+=(1,0)$ and $e_-=(0,1)$. The scattering matrix
$S(\l)$ for the pair $H$ and $H_0=-iJ {d\/ dx}$ has the following
form
\[
\lb{Sm} S(\l)={1\/a}\ma 1& -  \ol b\\ b& 1\am (\l ),\qq \l\in \R,
\]
see e.e. \cite{DEGM82}. Here ${1\/a}$ is the transmission
coefficient and $-{\ol b\/a}$ (or  ${ b\/a}$) is the right (left)
reflection coefficient. We have
\[
\lb{5.11} a(\l )=\det (\p^+,\vp^-)=\p_1^{+}(0,\l ),\qqq \qqq b(\l
)=-\p_1^{-}(0,\l ).
\]
The function $a(\l)$ is analytic in the upper half-plane $\C_+$ and
has an analytic extension in the whole complex plane $\C$.  All
zeros of $a$ lie in $\C_-$ (on the so-caleed non-physical sheet). We
denote by $(\l_n)_1^{\iy}$ the sequence of zeros of $a$
(multiplicities counted by repetition), so arranged that
$0<|\l_1|\leq |\l_2|\leq |\l_2|\leq \dots$.
 By the definition, the zero $\l_n\in \C_-$ of $a$ is a resonance.
  The multiplicity of the resonance is the multiplicity of
  the corresponding zero of $a$.

We recall some facts from \cite{IK12}. Let $R_0(\l)=(H_0-\l)^{-1}$
and $\gF(\l)=V_1R_0(\l)|q|^{1\/2}$ and $V=V_1 |q|^{1\/2}$. Note that
$\gF(\l), \Im \l\ne 0$ is the Hilbert-Schmidt operator, but the
operator $\gF(\l)$ is not trace class. In this case we define the
modified Fredholm determinant $D(\l)$ by
\[
D(\l)=\det\left[ (I+\gF(\l)) e^{-\gF(\l)}\right], \qqq \l\in \C_+.
\]
  We formulate now some results about resonances from \cite{IK12}:
{\it The determinant $D(\l)$ is analytic in $\C_+$ and has an
analytic extension into the whole complex plane $\C$ and
 $D=a$. Thus all zeros of $D$ are zeros of $a$, lie in $\C_-$ and
 satisfy
\[
\lb{ik} \cN(r,D)= {2r\/ \pi }(\g+o(1))\qqq as \qqq r\to\iy,
\]
where $\g>0$ is a diameter of the support of the potential $q$. Here
we denote the number of zeros of function $f$ having modulus  $\leq
r$ by $\cN (r,f)$, each zero being counted according to its
multiplicity. }
We formulate our main result.

\begin{theorem}
\lb{T1} Let the potential $q\in L^1(\R)$ and let $\supp q\ss
[0,\g]$,
 but in no smaller interval.  Then for each $p>1$ the following estimate hold true:
\[
\lb{r} \sum_{\Im \l_n< 0} {1\/|\l_n-i|^{p}}\le   {C Y_p\/\log
2}\rt({4\g\/\pi}+\int_\R |q(x)|dx\rt),
\]
where $C\le 2^5$ is an absolute  constant and  $Y_p=\int_\R
(1+x^2)^{-{p\/2}}dx, p>1$.
\end{theorem}

{\bf Remark.} 1) If  $\g+\int_\R |q(x)|dx\to 0$, then all resonances
go infinity.

2) The proof is  based on analysis of the function $a$ and
 the Carleson measure arguments \cite{C58}, \cite{C62}.
 We use harmonic  analysis  and Carleson's
Theorem (see  Theorems 1.56 and 2.3.9 in \cite{G81} and references
therein) about the Carleson measure. In fact, we use the approach
from \cite{K12}, where the estimates of resonances in terms of the
norm of potentials for 1D Schr\"odinger operators were obtained.
Note that in the case of the Dirac operator we obtain the sharper
estimate \er{r}.

3)   $C$ is an absolute  constant from Carleson's Theorem (
\cite{C58}, \cite{C62},   see Theorems 1.56 and   2.3.9,
\cite{G81}), see also \er{Ce}.

4) In fact, the estimate \er{r} gives a new global property of
resonance stability.

5)  The function $Y$ is strongly monotonic  on $(1,\iy)$  and
satisfies
\[
\lb{Yp} Y_2=\pi ,\qqq
Y_p=\ca {\sqrt{2\pi}\/\sqrt p}(1+o(1))  & as \ p\to \iy\\
             {2+o(1)\/p-1} & as \ p\to  1\ac.
\]
These properties of the function $Y_p$ is discussed in \cite{K12}.
In particular, asymptotics \er{Yp} are proved. Thus we can control
the RHS of \er{r} at $p\to 1$ and large $p\to \iy$. Note we take
$p>1$, since  the asymptotics \er{ik} implies the simple fact
$\sum_{n\ge 1} {1\/|\l_n|}=\iy$, see p. 17 in \cite{L93}.

\medskip

Resonances for the multidimensional case were studied by Melrose,
Sj\"ostrand, and Zworski and other, see \cite{M83}, \cite{Z89},
\cite{SZ91}) and references therein. We discuss the one dimensional
case. A lot of papers are devoted to the resonances for the 1D
Schr\"odinger operator, see Froese \cite{F97}, Simon \cite{S00},
Zworski \cite{Z87}, \cite{K11} and references therein. We recall
that Zworski \cite{Z87} obtained the first results about the
asymptotic distribution of resonances for the Schr\"odinger operator
with compactly supported potentials on the real line. Different
properties of resonances were determined in \cite{H99}, \cite{S00},
\cite{Z87} and  \cite{K04}, \cite{K05}, \cite{K11}. Inverse problems
(characterization, recovering, plus uniqueness) in terms of
resonances were solved by Korotyaev  for the Schr\"odinger operator
with a compactly supported potential on the real line \cite{K05} and
the half-line \cite{K04},   see also \cite{Z02} about uniqueness.

The "local resonance" stability problem was considered in
\cite{K04s}. It was proved that:  if $\vk=(\vk)_1^\iy$ is a sequence
of eigenvalues and resonances of the Schr\"odinger operator with
some compactly supported potential $q$ on the half-line and $\sum
_{n\ge 1}n^{2\ve}|\wt\vk_n-\vk_n|^2<\iy$ for some  sequence
$\wt\vk=(\wt\vk_n)_1^\iy$ and $\ve>1$, then $\wt\vk$ is a sequence
of eigenvalues and resonances of a Schr\"odinger operator with for
some unique real compactly supported potential $\wt q$. Another type
of the local resonance stability problem  was studied in
\cite{MSW10}.

Consider the Schr\"odinger operator $\cH = -\D - V$ acting in
$L^2(\R^d), d\ge 1$, where the potential $V\ge 0$ decreases
sufficiently fastly at infinity. The negative part of the spectrum
of $\cH$ is discrete and let $E_n<0, n\ge 1$ be the corresponding
increasing sequence of eigenvalues, each eigenvalue is counted
according to its multiplicity. This sequence is either finite or
tends to zero. Lieb and Thirring \cite{LT} proved inequalities of
the type
$$
\sum_{n\ge 1}|E_n|^\t\le C_{\t,d}\int_{\R^d}V^{{d\/2}+\t}dx, \qqq
$$
for some positive $\t$. There are a lot of papers about such
inequalities, see \cite{LS10}, \cite{LW00} and references therein.
In fact \er{r} is the Lieb-Thirring type inequalities for resonances of
the Dirac operator.

\

\section{Proof}
\setcounter{equation}{0}

\

{\bf 2.1. Estimates for entire functions.} An entire function $f(z)$
is said to be of exponential type if there is a constant $\a$ such
that $|f(z)|\leq\const e^{\a|z|}$ everywhere.  The infimum of the
set of $\a$ for which such inequality holds is called the type of
$f$.

\no {\bf Definition.} {\it Let $\cE_\g, \g>0$ denote the space of
exponential type  functions $f$,  which satisfy
\[
\lb{CE}
|f(\l)|\le e^{A+\g(|\Im \l|-\Im \l)} \qqq \qq  \forall \ \l\in \C,
\]
\[
\lb{CE1}
|f(\l)|\ge 1 \qqq \forall \ \ \l\in \R,
\]
for some constants $A=A(f)\ge 0$.}

In the proof of Theorem \ref{T1} we need some properties of the zeros of
$f\in \cE_\g$ in terms of the Carleson  measure. Recall that {\it a
positive Borel measure $M$ defined in $\C_-$ is called the
Carleson measure if there is
 a constant $C_M$ such that for all $(r,t)\in \R_+\ts\R$
\[
   \lb{1.31}
   M(D_-(t,r))\leq C_Mr,\ \ \ {\rm where}\ \ \
   \ D_-(t,r)\ev\{z\in \C_-: |z-t|<r \},
\]
here $C_M$ is the Carleson  constant independent of $(t,r)$.}

For an entire function $f$ with zeroes $\l_n, n\ge 1$  we define an associated measure by
\[
\lb{MO}
 d\O (\l,f)=\sum_{\Im \l_n\le 0} \d (\l-\l_n+i)d\m d\e,\ \ \ \ \l=\m+i\e\in \C_-.
\]

In order to prove Theorem \ref{T1} we need following results.

  \begin{theorem}
\lb{TC}
 Let $f\in \cE_\g, \g>0$. Then

 i) for each $r>0$ the following
estimate hold true:
\[
\lb{L1} \cN (r,f)\leq {1\/\log 2}\rt({4r\g\/\pi}+ A\rt).
\]
ii) $d\O (\l,f)$ is the Carleson measure and satisfies
\[
\lb{L2} \O (D_-(t,r),f)\le\cN (r,f(t+\cdot))\le {r\/\log
2}\rt({4\g\/\pi}+ A\rt) \qqq \forall \ (r,t)\in \R_+\ts \R.
\]
iii) For each $p>1$ the following estimates hold true:
\[
\lb{L3} \sum_{n\ge 1} {1\/|\l_n-i|^{p}}\le   {C Y_p\/\log
2}\rt({4\g\/\pi}+A\rt),
\]
where  $C\le 2^5$ is an absolute  constant and  $Y_p=\int_\R
{dx\/(1+x^2)^{p\/2}}, p>1$.
\end{theorem}

{\bf Proof.} i) Let $F=e^{-i\g\l}f, \l=re^{i\f}$. Then the Jensen
formula implies (see 2 p. in \cite{Koo88})
\[
\lb{X1}
\log |f(0)|+\int _0^r{\cN (t,f)\/ t}dt= {1\/ 2\pi }\int
_0^{2\pi}\log |F(re^{i\f})|d\f,
\]
since $\cN (t,f)=\cN (t,F)$. Using the estimate \er{CE} we obtain
$$
\log |F(re^{i\f})|\le \g r |\sin \f|+ A, \ \ \  \l=re^{i\f},\ \ \ \f
\in [0, 2\pi ],
$$
which yields
\[
\lb{X2} {1\/ 2\pi }\int _0^{2\pi}\log |F(re^{i\f})|d\f \le {2 \g\/
\pi}r+ A.
\]
Substituting the estimate \er{X2} into the identity \er{X1}
together  with the simple estimate
$$
\int _0^r{\cN (t)dt\/
t}\geq \cN \rt({r\/2}\rt)\int _{r/2}^r{dt\/ t}=\cN \rt({r\/2}\rt)\log 2,
$$
we obtain \er{L1}, since $|f(0)|\ge 1$.

ii) Let $r\le 1, t\in \R$. Then by the construction of $\O
(\cdot,f)$,  we obtain $\O (D_-(t,r),f) =0$.

Let  $r>1, t\in \R$. Then due to \er{L1}, the measure $\O (\cdot,f)$ satisfies
$$
\O (D_+(t,r),f)\le\cN (r,f(t+\cdot))\le {1\/\log 2}\rt({4r\g\/\pi}+
A\rt)\le {r\/\log 2}\rt({4\g\/\pi}+ A\rt).
$$
Thus $\O (\cdot,f)$ is the Carleson measure with the Carleson
$C_\O={1\/\log 2}\rt({4\g\/\pi}+ A\rt)$.

iii)  In order to show \er{L3} we recall the Carleson result (see p.
63, Theorem 3.9, \cite{G81}):

Let $f$ be analytic on $\C_-$. For $0<p<\iy$  we say $f\in
\mH_p=\mH_p(\C_-)$ if
$$
\sup_{y<0}\int_\R|f(x+iy)|^pdx=\|f\|_{\mH_p}<\iy
$$
Note that the definition of the Hardy space $\mH_p$ involve all
$y>0$, instead of small only value of $y$, like say, $y\in (0,1)$.
We define the Hardy space $\mH_p$ for the case $\C_-$, since below
we work with functions on $\C_-$.

{\it If $M$ is a Carleson measure, then the following estimate
holds:
\[
\lb{Ce}
  \int_{\C_-} |f|^pdM\leq C C_M \|f\|_{\mH_p}^p\qqq
  \forall \qq f\in \mH_p, \ p\in (0,\iy ),
\]
where  $C_M$ is the so-called Carleson constant in \er{1.31} and
$C\le 2^5$ is an absolute constant.}

 In order to prove \er{L3} we  take the functions $f(\l)={1\/\l-i}$.
Estimate \er{Ce} yields
\[
\lb{3.19}
 \int_{\C_-} |f(\l)|^pdM=\sum_{n\ge 1} {1\/|\l_n-i|^p} \leq
C C_M \|f\|_{\mH_p}^p, \qq p\in (1,\iy).
\]
where we have the simple identity
\[
\lb{3.19a} \|f\|_{\mH_p}^p=\int_\R {dt\/|t-i|^p}=\int_\R
{dt\/(t^2+1)^{p\/2}}=Y_p.
\]
 Combine \er{3.19} and \er{3.19a} we obtain \er{L3}.
\hfill  \BBox

{\bf 2.1. Estimates of resonances.} We consider the function $a$
given by \er{5.11} under the condition that the potential $q$
satisfies $\supp q\ss [0,\g]$. The solution $\p^+=(\p_1^+,\p_2^+)$
satisfy the integral equations
\[
\p^{+}(x,\l )=e^{i\l x}e_{1}+\int _x^\g iJe^{i\l (x-t)J} V(t) \p^{+}(t,\l )dt,
\]
where
\begin{equation}
\lb{5.6} iJe^{i\l (x-t)J} V(t) = i\ma 0& q(t) e^{i\l (x-t)}\\
-\ol{q}(t)e^{-i\l (x-t)}&0\am .
\end{equation}
For $\c=e^{-i\l x}\p_{1}^+(x,\l )$ using (\ref{5.6}) we obtain
$$
 \c(x,\l)=1+i\int _x^\g e^{-i\l s}q(s) \p_{2}^+(s,\l )ds,
$$
$$
 \p_{2}^+(s,\l )=-i\int _s^\g e^{i\l (2t-s)}\ol{q}(t) \c(t,\l )dt.
$$
Then $\c(x,\l )$ satisfies the following integral equation
\[
\lb{x1}
\begin{aligned}
 \c(x,\l )=1+ \int _x^\g q(t_1)dt_1 \int _{t_1}^\g e^{i2\l (t_2-t_1)}
 \ol{q}(t_2)\c(t_2,\l)dt_2.
\end{aligned}
\]
We have the standard formal iterations given by
\[
\lb{ien}
\c (x,\l )=1+\sum _{n\geq 1}\c_n(x,\l ),\ \ \ \
\c_n(x,\l )=\int _x^\g q(t_1)dt_1 \int _{t_1}^\g e^{i2\l (t_2-t_1)}
\ol{q}(t_2)\c_{n-1}(t_2,\l)dt_2,
\]
where $\c_0(\cdot ,\l )=1$. Due to \er{5.11} we get $a(\l )=\c (0,\l
)$,  which yields
\[
\lb{iea}
a(\l )=1+\sum _{n\geq 1}a_n(\l ),\qqq
  a_n(\l)=\c _n(0,\l).
\]
We will describe these iterations and the function $a$.

\begin{lemma}
\lb{Ta}
Let $q\in L^1(\R)$ and  $\supp q\ss [0,\g]$. Then the function
$a\in \cE_\g$ and satisfies
\[
\lb{an}
  |a_n(\l )|\leq e^{\g(|\e|-\e)}{\|q\|_1^{2n}\/ (2n)!},\ \
   \forall \ n\geq 1,
\]
\[
\lb{a}
  |a(\l)|\leq e^{\g(|\e|-\e)}\ch \|q\|_1,
\]
\[
\lb{a1}
|a(\l)-1|\leq e^{\g(|\e|-\e)}(\ch \| q\|-1).
\]
 where $\e=\Im\l$ and   $\|q\|_1=\int_\R |q(x)|dx.$

\end{lemma}
\no {\bf Proof.}  Let $\e_-={(|\e|-\e)\/2}$. Then using \er{ien} we
obtain
$$
a_n(\l)=\int\limits_{0=t_0<t_1< t_2<...< t_{2n}}
\lt(\prod\limits_{1\le j\le n}  q(t_{2j-1})\ol{q}(t_{2j}) e^{i2\l
(t_{2j}-t_{2j-1})}
\rt)dt_1dt_2...dt_{2n},
$$
which yields
\[
\begin{aligned}
\lb{y4} |a_n(\l)|\le \int\limits_{0=t_0<t_1< t_2<...< t_{2n}}
\lt(\prod\limits_{1\le j\le n}e^{2\e_- (t_{2j}-t_{2j-1})}
|q(t_{2j-1})q(t_{2j})|\rt)dt_1dt_2...dt_{2n}\\
\le \int\limits_{0<t_1< t_2<...< t_{2n}} \lt(\prod\limits_{1\le j\le
{2n}}
|q(t_j)|\rt) e^{2\e_- t_{2n}}dt_1dt_2...dt_{2n}\\
\le e^{\g (|\e|-\e)}\int\limits_{0<t_1< t_2<...< t_{2n}}
|q(t_1)q(t_2)....q(t_{2n})| dt_1dt_2...dt_{2n}\le e^{\g (|\e|-\e)}
{\|q\|^{2n}\/(2n)!},
\end{aligned}
\]
which yields \er{an}.

This shows that the series \er{iea} converges uniformly on bounded
subsets of $\C$. Each term of this series is an entire function.
Hence the sum is an entire function. Summing the majorants we obtain
estimates  \er{a} and \er{a1}. Thus the function  $a$ is entire.

The S-matrix is unitary, see \cite{DEGM82}, then from \er{Sm} we
have the well-known fact $|a(\l)|^2=|b(\l)|^2+1\ge 1$ for all $\l\in
\R$.  Then due to \er{a1} we deduce that $a$ belongs to $\cE_\g$.
\hfill \BBox

{\bf Proof of Theorem \ref{T1}.} Recall  that by Lemma \ref{Ta}, the
function $a\in \cE_\g$ with $A=\|q\|_1$, since $a$  satisfies $
|a(\l)|\le e^{\g(|\e|-\e)+\|q\|_1}$, for all $\l\in\C$ with $ \e=\Im
\l$. Thus estimate \er{L3}  gives the main estimate \er{r} in
Theorem \ref{T1}.
 \hfill $\BBox$

\


\no  {\bf Acknowledgments.}\small
 Various parts of this paper were written during Korotyaev's stay in
 the Mathematical Institute of University of Aarhus.
  He is grateful to the institute for the hospitality.
  He is also grateful to Iliya Vidensky and especially to Alexei
Alexandrov  (St. Petersburg)  for stimulating discussions about the
Carleson Theorem.

This work was supported by the
Ministry of education and science of the Russian Federation, state
contract 14.740.11.0581 and the RFFI grant "Spectral and asymptotic
methods for studying of the differential operators" No 11-01-00458.


\end{document}